\documentclass[12pt,notitlepage]{article}

\usepackage{chngcntr}
\usepackage[margin=0.9in]{geometry}
\usepackage[parfill]{parskip}   
\usepackage[nodisplayskipstretch]{setspace}
\usepackage{enumerate}
\usepackage{bm}
\usepackage{amssymb}
\usepackage{amsmath}

\usepackage{amsthm}
\usepackage{graphicx}
\usepackage{subfigure}
\usepackage{makeidx}
\usepackage{multicol}
\usepackage{comment}
\usepackage{cite}
\usepackage{color}
\usepackage{chngcntr}
\definecolor{Gray}{gray}{0.9}
\newtheorem{theorem}{Theorem}
\newtheorem*{mproof}{Proof}

\newtheorem{proposition}{Proposition}
\newtheorem{corollary}{Corollary}
\newtheorem{remark}{Remark}

\makeatletter
\renewenvironment{thebibliography}[1]{%
  \@xp\section\@xp*\@xp{\refname}%
 \normalfont\footnotesize\labelsep .5em\relax
  \renewcommand\theenumiv{\arabic{enumiv}}\let\p@enumiv\@empty
  \vspace*{-20pt}
  \list{\@biblabel{\theenumiv}}{\settowidth\labelwidth{\@biblabel{#1}}%
    \leftmargin\labelwidth \advance\leftmargin\labelsep
    \usecounter{enumiv}}%
  \sloppy \clubpenalty\@M \widowpenalty\clubpenalty
  \sfcode`\.=\@m
}{%
  \def\@noitemerr{\@latex@warning{Empty `thebibliography' environment}}%
  \endlist
}
\makeatother

\date{}

\begin{document}

\title{\large\bf 
$\,$\\ 
\vspace{-34pt}
On The Number Of Labeled Bipartite Graphs}

\vspace{-44pt}
\author{\ Abdullah Atmaca\footnote{Department of Computer Engineering, Bilkent University, Ankara, Turkey.}   \,
  and A. Yavuz Oru\c{c}\footnote{Department of Electrical and Computer Engineering, University of Maryland, College Park, MD 20742.}
}

\maketitle
\nopagebreak[4]
\setcounter{page}{1}

\vspace{-20pt}
\begin{center}
{\bf Abstract}
\end{center}

\vspace{-1pt}\noindent
Let  $I$ and $O$ denote two sets of vertices, where $I\cap O =\O$,  $|I| = n$, $|O| = r$,  and  $B_u(n,r)$ denote the set of unlabeled graphs whose edges connect vertices in $I$ and $O$. Recently, it was established in\cite{atmacaoruc2018} that the following two-sided equality holds,
\vspace{-1pt}
\begin{equation}
\label{mainResult}
\displaystyle  \frac{\binom{r+2^{n}-1}{r}}{n!}\, \le\, |B_u(n,r)|\, \le\,  2\frac{\binom{r+2^{n}-1}{r}}{n!},\, n <  r.\nonumber
\end{equation}

\vspace{-1pt}\noindent
and exact formulas were provided in~\cite{atmaca2017size}  for $|B_u(2,r)|$ and $|B_u(3,r)|.$ 
In this paper,  these results are extended to various families of labeled bipartite graphs.

\vspace{5pt}\noindent

\vspace{5pt}\noindent

\vspace{-10pt}

\section{Introduction}

\begin{figure}[t!]
\centering{
\vspace{-5pt}
\includegraphics[scale=1.2]{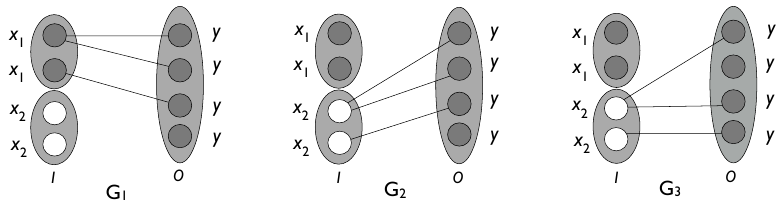}
\vspace{-5pt}
\caption{Examples of equivalent and non-equivalent left-set-labeled $(4,4)$-bipartite graphs.}
\label{figure1}
}
\end{figure}

\vspace{-10pt}\noindent
Exact and asymptotic enumerations of graphs that satisfy certain properties have been a key problem in graph theory~\cite{harary1958number,harary1963enumeration,wright1971graphs,wright1972number, harrison1973number,wright1974graphs,wright1974asymmetric,wright1974two,wright1976proportion,hanlon1979enumeration,bollobas1982asymptotic,gainer2014enumeration,atmaca2017size,atmacaoruc2018}. A number of results have been reported on the problem of counting graphs in which number of vertices and/or edges are fixed\cite{wright1971graphs,wright1972number,wright1974graphs,wright1974asymmetric,wright1974two,bollobas1982asymptotic}. One such particular enumeration deals with bipartite graphs with $n$ left vertices and $r$ right vertices in which left (right) vertices are unlabeled, i.e., they are indistinguishable\footnote{The terms left (right) vertices are used only for notational convenience, and otherwise have no bearing on the validity of the results presented in the paper.}. Using Polya's counting theorem, it was recently shown in~\cite{atmacaoruc2018} that the number of such graphs is given by $c\frac{\binom{r+2^{n}-1}{r}}{n!}$ if  $n < r,$ and by $c\frac{\binom{n+2^{r}-1}{n}}{r!}$ if $r < n,$ where $1\le c \le 2.$ In addition, exact formulas were given for the number of such bipartite graphs with two and three left vertices  and $r$ right vertices, $r\ge 3,$ in~\cite{atmaca2017size}.  Unlabeled graphs constitute the smallest set of distinct bipartite graphs as any two such graphs are considered equivalent under the most unconstrained form of equivalence, where only the adjacencies of vertices are enforced without distinguishing between the vertices within $I$ and those within $O.$ For example, all three  bipartite graphs shown in Figure~\ref{figure1} are equivalent under the usual definition of graph isomorphism. However, only  $G_2$ and $G_3$ are equivalent and $G_1$ is not equivalent to either  $G_2$ or $G_3$, assuming that subsets of left vertices are distinguished from each other, while right vertices are indistinguishable in this example. Differentiating between bipartite graphs when left and/or right vertices represent different objects as subsets is useful in some applications, for example in classifying and enumerating distinct calls (communications) between a set of callers and a set of receivers, based on the identities of callers and receivers in a combinatorial call model and classifying corresponding switching networks into various families such as concentrators, superconcentrators and generalizers~\cite{oruc2016}. In this paper, we provide formulas for the cardinalities of various families of such labeled bipartite graphs.

The rest of the paper is organized as follows. The next section formalizes the notion of labeled bipartite graphs and counting problems for such graphs.  We provide exact formulas for left-set-label $(n,r)$-bipartite graphs in Section~\ref{leftSetLabel}, and derive  similar formulas for set-label $(n,r)$-bipartite graphs in Section~\ref{setLabel},  by fixing one of $n$ and $r$ to 2 and 3. Sections 5 and 6 present lower and upper bounds for the number of left and set-labeled $(n,r)$-bipartite graphs  for all $n, r \ge 1.$ 

$\,$\vspace{-10pt}
\section{Preliminaries}
\label{preliminaries}

\vspace{-5pt}\noindent
 Let $ G = (I\, \oplus\,  O, E)$  be a bipartite graph with two non-overlapping sets of (left) vertices, $I$ and (right) vertices, $O,$ and a set of edges, $E$ whose elements join vertices in $I$ with those in $O.$ If $|I| = n$ and $|O| = r$, $G$ will be called an  $(n,r)$-bipartite graph. Labeling left (right) vertices in an $(n,r)$-bipartite graph partitions them into subsets so that any two left (right) vertices belong to the same subset if and only if they have the same label. Two special cases occur when all left (right) vertices have the same label or when they all have different labels. The first case implies that the left (right) vertices in $G$ are entirely indistinguishable, and the second case implies that the left (right) vertices in $G$ are entirely distinguishable.  For example, in all three bipartite graphs in Figure~\ref{figure1}, the top two left vertices  are distinguishable from the bottom two left vertices, and none of the right vertices  is distinguishable from the other right vertices.
 
 An $(n,r)$-bipartite graph $G$ is called {\em unlabeled} if both of its left and right vertices are entirely indistinguishable.  It is called  {\em left-set-labeled} if its left vertices are distinguishable up to the subsets to which they belong and right vertices are entirely indistinguishable. It is called {\em right-set-labeled} if its right vertices are distinguishable up to the subsets to which they belong and its left vertices are entirely indistinguishable.  Furthermore, it is called {\em set-labeled} if both of its left  and right vertices are distinguishable up to subsets to which they belong.  All three graphs in Figure~\ref{figure1} are left-set-labeled graphs.
The sets of these four types of labeled bipartite graphs will be denoted by $ B_u(n,r), B_x(n,r), B_y(n,r),$ and $B_{x,y}(n,r),$ in that order. 
 
The subset of all left vertices in a bipartite graph, $G$ with non-zero degrees is called its  {\em connected domain} and denoted by $I_c(G),$ and the subset of right vertices with non-zero degrees is called the {\em connected co-domain} of $G,$ and denoted by $O_c(G).$ Two  bipartite graphs $G_1 = (I\oplus O, E_1)$ and  $G_2 = (I\oplus O, E_2)$ are  {\em left-set-labeled equivalent} if there exist bijections $\alpha: I\rightarrow I$ and $\beta: O\rightarrow O$ such that: (i) for any $x\in I$ and $y\in O,$  $(x,y)\in E_1$ if and only if $(\alpha(x),\beta(y))\in E_1$ and  (ii) $I_c (G_1) =  I_c(G_2).$ Similarly, $G_1$ and $G_2$ are  {\em right-set-labeled equivalent}  if there exist bijections $\alpha: I\rightarrow I$ and $\beta: O\rightarrow O$  such that (i) holds and  $O_c (G_1) =  O_c(G_2),$  and they are set-labeled equivalent if there exist bijections $\alpha: I\rightarrow I$ and $\beta: O\rightarrow O$  such that (i) holds, and both  $I_c (G_1) =  I_c(G_2)$ and  $O_c (G_1) =  O_c(G_2).$ The bipartite graphs, $G_2$ and $G_3$  in Figure~\ref{figure1} are left-set-labeled equivalent as they satisfy both conditions, but neither $G_1$ and $G_2$ nor $G_1$ and $G_3$ are, as they do not satisfy part (ii) of the definition. In addition,  $G_1$ and $G_2$ are right-set-labeled equivalent, but they are not set-labeled equivalent.  

It is not difficult to verify that the following remark holds.
 \begin{remark}
 {\rm $\,$\\
 
 \vspace{-12pt}
(i) $\,\,\,B_u(n,r)\subseteq B_x(n,r),$ 

(ii)  $\,B_u(n,r)\subseteq B_y(n,r),$ 

(iii)  $B_x(n,r), B_y(n,r) \subseteq  B_{x,y}(n,r).$ \qed
 }
 \end{remark}

\vspace{5pt}\noindent
Moreover, the following result was established in~\cite{atmaca2017size} for unlabeled bipartite graphs.

\vspace{1pt}
\begin{theorem}
{\rm \label{theorem2r}

$\,$\\

\vspace{-50pt}
\begin{eqnarray}
\label{TheOneCase}
\text{(a)}\, |B_u(1,r)| &=& r+1,\\\\
\label{TheTwoCase}
\text{(b)}\, |B_u(2,r)|&=& \frac{2r^{3}+15r^{2} + 34r + 22.5 + 1.5\left ( -1 \right )^{r}}{24},\\\\
\label{TheThreeCase}
\text{(c)}\, |B_{u}\left( 3,r \right)| &=& 
\left\{\begin{matrix}
\frac{1}{6}\left [ A(r) + \frac{2(r^{3}+12r^{2}+45r+54)}{54} \right ] &\!\!\!\!\!\! \text{ if }\, r  \bmod\!
\text{ } 3 = 0, \\
\\
\frac{1}{6}\left [A(r) + \frac{2(r^{3}+12r^{2}+45r+50)}{54} \right ] &\!\!\!\!\!\! \text{ if }\, r  \bmod\!\!
\text{ } 3 = 1, \\
\\
\frac{1}{6}\left [ A(r) + \frac{2(r^{3}+12r^{2}+39r+28)}{54} \right ] &\!\!\!\!\!\! \text{ if }\, r  \bmod\!\!
\text{ } 3 = 2, \!\!\!\! \end{matrix}\right.,
\end{eqnarray}
 
\vspace{4pt}\noindent 
where 
$A(r) =  \binom{r+7}{r} + \frac{3\left ( r+4 \right )\left ( 2r^{4}+32r^{3}+172r^{2} + 352r
	+ 15\left ( -1 \right )^{r} +225 \right )}{960},$  and
	\vspace{-1pt}
		
\begin{equation}
\!\!\!\!	\!\!\!\!	\!\!\!\!	\!\!\!\!	\!\!\!\!	\!\!\!\!	\!\!\!\!	\!\!\!\!	\!\!\!\!	\!\!\!\!	\text{(d)}\,	\frac{\binom{r+2^{n}-1}{r}}{n!}\, \le\, |B_u(n,r)|\, \le\,  2\frac{\binom{r+2^{n}-1}{r}}{n!},\, n <  r. \label{generalCase}\qed
	\end{equation}
}
\end{theorem}

\vspace{-1pt}\noindent
In what follows, we count the number of bipartite graphs in $B_x(n,r)$  when (i) $n = 2, n = 3,$ and $r\ge 1,$ (ii)  $r = 2, r = 3,$ and $n\ge 1,$  and provide lower and upper bounds for all $n,r\ge 1.$ We also count the number of bipartite graphs in $B_{x,y}(n,r)$ when $r = 2, r = 3$ and $n \ge 1.$   Right-set-labeled bipartite graphs in  $B_y(n,r)$ are counted similarly and their counting is omitted.

\vspace{-4pt}
\section{Left-Set-Labeled Bipartite Graphs}
\label{leftSetLabel}

\vspace{-7pt}\noindent
Let $B_x(n,r,i)$ denote the set of $(n,r)$-left-set-labeled bipartite graphs in each of which the degrees of exactly $i$ left vertices are greater than 0. Let $\overline{B}_x(i,r)$ denote the set of all $(i,r)$-unlabeled bipartite graphs in each of which the degrees of all left vertices are greater than 0. Given that $i$ left vertices can be fixed in any  one of $\binom{n}{i}$ ways out of $n$ left vertices, and placing each bipartite graph,  some fixed  $i$ left vertices of which all have non-zero degrees  in $B_x(n,r,i)$  in one-to-one correspondence with  a bipartite graph in $\overline{B}_x(i,r),$   we have 
\vspace{-10pt}
\begin{eqnarray}
	|B_x(n,r,i)| &=& \binom{n}{i} |\overline{B}_x(i,r)|,\\
	|B_x(n,r)| &=& 1 + \sum_{i=1}^n |B_x(n,r,i)|,\\
	|B_x(n,r)| &=& 1 + \sum_{i=1}^n \binom{n}{i} |\overline{B}_x(i,r)|. \label{bxnrFormula}	
\end{eqnarray}

\vspace{-10pt}\noindent
By Theorem~\ref{theorem2r},  $|B_u(1,r)|=r+1$. Dropping the  bipartite graph with no edges yields $|\overline{B}_x(1,r)| = r$, and hence by Eqn.~\ref{bxnrFormula}, $|B_x(1,r)| = r+1$. These $r+1$ graphs in  $B_x(1,r)$ all have a single left vertex and can have $0, 1,\cdots r$ right vertices as its neighbors.

To compute $|B_x(n,r)|$, $n \ge 2$ we note the following identities.
\vspace{-10pt}
\begin{eqnarray}
|B_u(n,r)| = 1 + \sum_{i=1}^n |\overline{B}_x(i,r)|,  \\
|B_u(n-1,r)| = 1 + \sum_{i=1}^{n-1} |\overline{B}_x(i,r)|.
\end{eqnarray}

\vspace{-15pt}\noindent
Hence
\begin{equation}
|\overline{B}_x(n,r)| = |B_u(n,r)| - |B_u(n-1,r)|.
\label{bxnrRecFormula}
\end{equation}

\noindent 
Now we can calculate $|\overline{B}_x(2,r)|$ using Eqn.~\ref{bxnrRecFormula} and Theorem~\ref{theorem2r}(Eqn.~\ref{TheTwoCase}),
\vspace{-1pt}
\begin{eqnarray}
|\overline{B}_x(2,r)| & = & |B_u(2,r)| - |B_u(1,r)|, \\
 &= & \frac{2r^{3}+15r^{2} + 34r + 22.5 + 1.5\left ( -1 \right )^{r}}{24} - (r+1),
\label{bbarx2r}
\end{eqnarray}

\vspace{-1pt}\noindent
and similarly $|\overline{B}_x(3,r)|$ using Eqn.~\ref{bxnrRecFormula} and Theorem~\ref{theorem2r}(Eqn.~\ref{TheThreeCase}),

\vspace{-2pt}\begin{eqnarray}
|\overline{B}_x(3,r)| &=& |B_u(3,r)| - |B_u(2,r)|,  \\ \\
&=& \left\{\begin{matrix}
\frac{1}{6}\Bigg[ \binom{r+7}{7} + \frac{3\left ( r+4 \right )\left ( 2r^{4}+32r^{3}+172r^{2} + 352r
	+ 15\left ( -1 \right )^{r} +225 \right )}{960} + \frac{2(r^{3}+12r^{2}+45r+54)}{54} \Bigg] \\  - \frac{2r^{3}+15r^{2} + 34r + 22.5 + 1.5\left ( -1 \right )^{r}}{24}  \text{ if }\, r  \bmod\!\!
\text{ } 3 = 0, \\
\\
\frac{1}{6}\Bigg[ \binom{r+7}{7} + \frac{3\left ( r+4 \right )\left ( 2r^{4}+32r^{3}+172r^{2} + 352r
	+ 15\left ( -1 \right )^{r} +225 \right )}{960} + \frac{2(r^{3}+12r^{2}+45r+50)}{54} \Bigg] \\ - \frac{2r^{3}+15r^{2} + 34r + 22.5 + 1.5\left ( -1 \right )^{r}}{24} \text{ if }\, r  \bmod\!\!
\text{ } 3 = 1, \\
\\
\frac{1}{6}\Bigg[ \binom{r+7}{7} + \frac{3\left ( r+4 \right )\left ( 2r^{4}+32r^{3}+172r^{2} + 352r
	+ 15\left ( -1 \right )^{r} +225 \right )}{960} + \frac{2(r^{3}+12r^{2}+39r+28)}{54}  \Bigg] \\  - \frac{2r^{3}+15r^{2} + 34r + 22.5 + 1.5\left ( -1 \right )^{r}}{24} \text{ if }\, r  \bmod\!\!
\text{ } 3 = 2. \!\!\!\! \end{matrix}\right. \nonumber \\
\text{ }
\label{bbarx3r}
\end{eqnarray}

\vspace{-5pt}\noindent
Using Eqns.~\ref{bxnrFormula} and \ref{bbarx2r}, we can easily calculate $|B_x(2,r)|$.

\begin{figure}[t!]
\centering{
\vspace{1pt}
\includegraphics[scale=0.87]{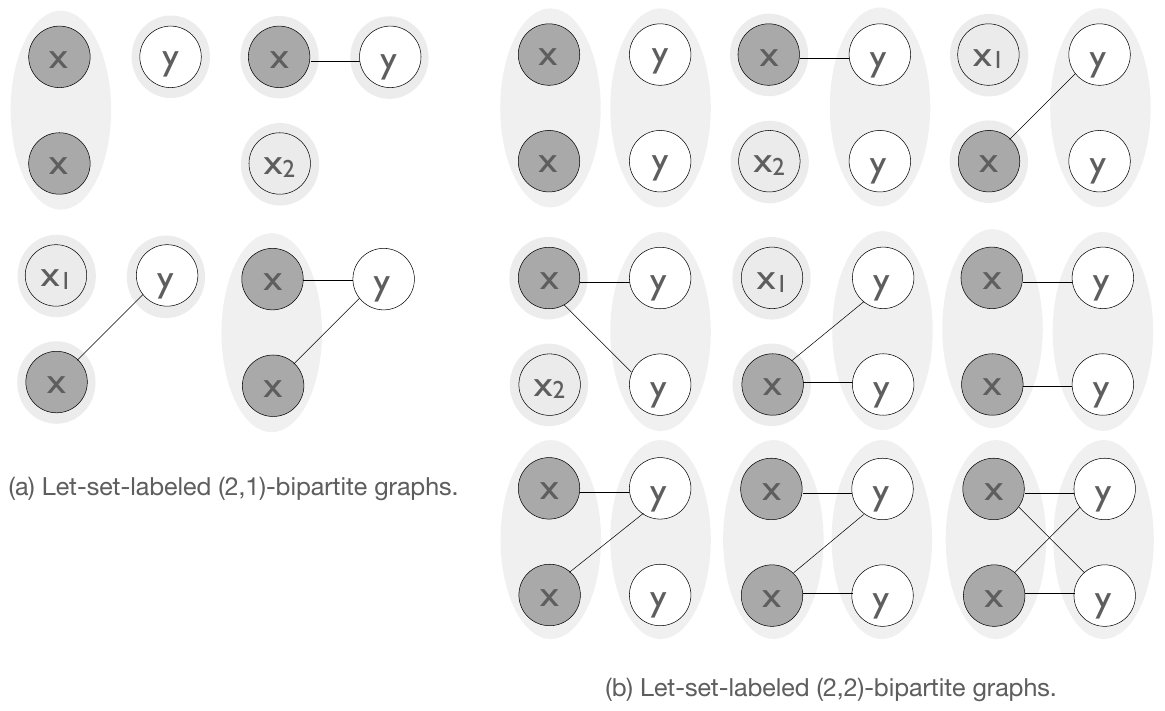}
\vspace{1pt}
\caption{Non-equivalent left-set-labeled (2,1) and (2,2)-bipartite graphs.}
\label{figure2}
}
\end{figure}

\begin{figure}[t!]
\centering{
\vspace{1pt}
\includegraphics[scale=0.28]{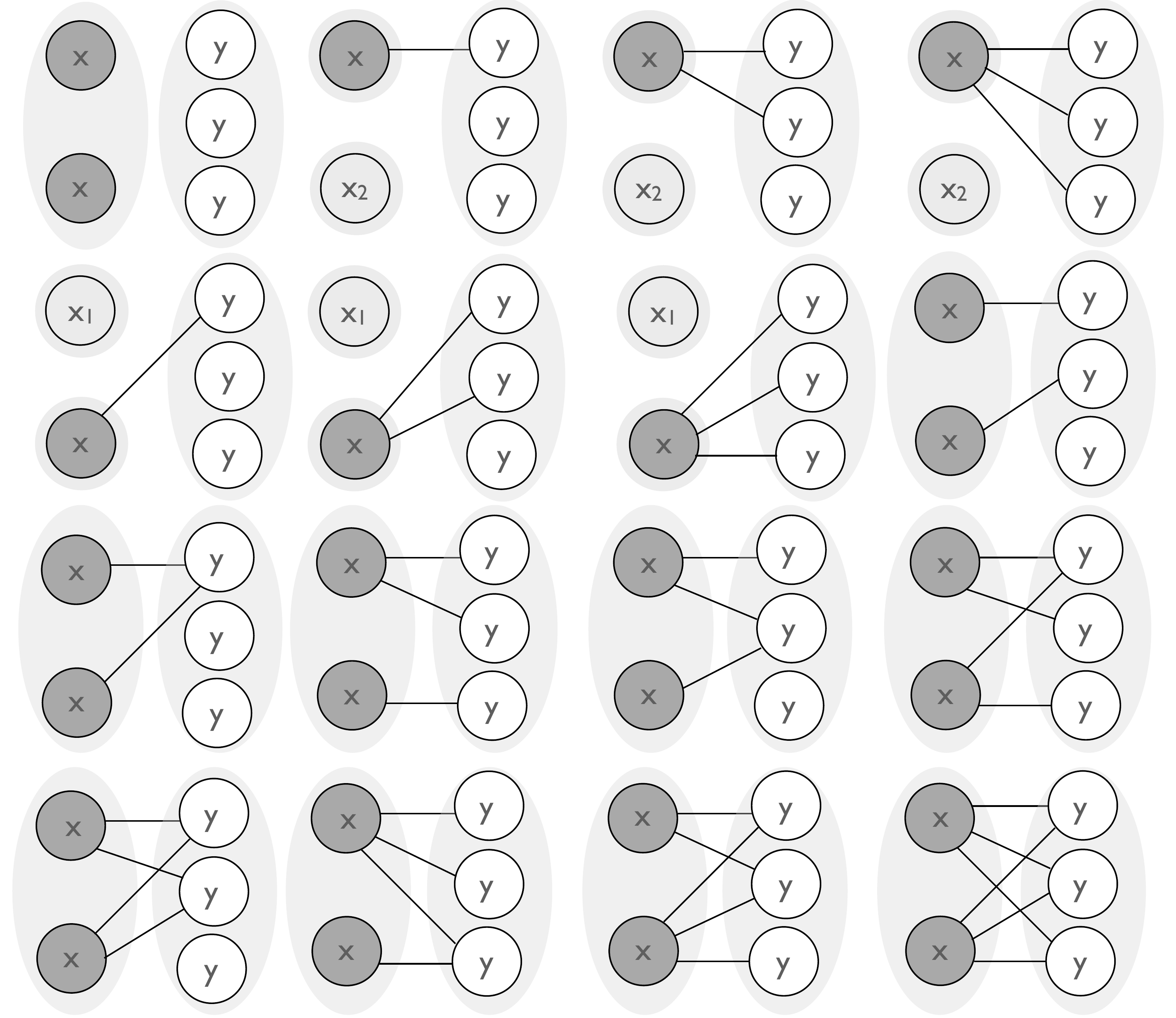}
\vspace{-4pt}
\caption{Non-equivalent left-set-labeled (2,3)-bipartite graphs.}
\label{figure3}
}
\end{figure}

\vspace{-5pt}
\begin{eqnarray}
|B_x(2,r)| &=& 1 + \sum_{i=1}^2 \binom{2}{i} |\overline{B}_x(i,r)|  = 1 + 2 |\overline{B}_x(1,r)| + |\overline{B}_x(2,r)|, \\
&=& 1 + 2r + \frac{2r^{3}+15r^{2} + 34r + 22.5 + 1.5\left ( -1 \right )^{r}}{24} - (r+1), \phantom{.......} \\
&=& \frac{2r^{3}+15r^{2} + 58r + 22.5 + 1.5\left ( -1 \right )^{r}}{24}.
\end{eqnarray}
\begin{remark}
{\rm
It follows that $|B_x(2,1)| = 4,$ $|B_x(2,2)| = 9$ and $|B_x(2,3)| = 16.$ The corresponding bipartite graphs in $B_x(2,1), B_x(2,2)$ and  $B_x(2,3)$ are shown in Figure~\ref{figure2} and Figure~\ref{figure3}. $||$
}
\end{remark}


\vspace{2pt}
Similarly, using Eqns.~\ref{bxnrFormula} and \ref{bbarx3r} we can calculate $|B_x(3,r)|$ as follows.
\vspace{-2pt}
\begin{eqnarray}
|B_x(3,r)| &\!\!\!\!=& \!\!\!\!\!  1 + \sum_{i=1}^3 \binom{3}{i} |\overline{B}_x(i,r)|,  \\
&\!\!\!\!=& \!\!\!\!\!  1 + 3 |\overline{B}_x(1,r)| + 3 |\overline{B}_x(2,r)| + |\overline{B}_x(3,r)|, \\
&\!\!\!\!=& \!\!\!\!\! 1 + 3r + 3 \Big( \frac{2r^{3}+15r^{2} + 34r + 22.5 + 1.5\left ( -1 \right )^{r}}{24} - (r+1) \Big) + |\overline{B}_x(3,r)|, \phantom{..........}
\end{eqnarray}

\vspace{-1pt}\noindent 
and replacing $|\overline{B}_x(3,r)|$ by the formula in Eqn.~\ref{bbarx3r} gives 

\vspace{1pt}
\begin{eqnarray}
|B_x(3,r)|=& \left\{\begin{matrix}
\frac{1}{6}\Bigg[ \binom{r+7}{7} + \frac{3\left ( r+4 \right )\left ( 2r^{4}+32r^{3}+172r^{2} + 352r
	+ 15\left ( -1 \right )^{r} +225 \right )}{960} + \frac{2(r^{3}+12r^{2}+45r+54)}{54} \Bigg] \\  + \frac{4r^{3}+30r^{2} + 68r -3 + 3 ( -1 )^{r}}{24}  \text{ if }\, r  \bmod\!\!
\text{ } 3 = 0, \\
\\
\frac{1}{6}\Bigg[ \binom{r+7}{7} + \frac{3\left ( r+4 \right )\left ( 2r^{4}+32r^{3}+172r^{2} + 352r
	+ 15\left ( -1 \right )^{r} +225 \right )}{960} + \frac{2(r^{3}+12r^{2}+45r+50)}{54} \Bigg] \\ + \frac{4r^{3}+30r^{2} + 68r -3 + 3 ( -1 )^{r}}{24} \text{ if }\, r  \bmod\!\!
\text{ } 3 = 1, \\
\\
\frac{1}{6}\Bigg[ \binom{r+7}{7} + \frac{3\left ( r+4 \right )\left ( 2r^{4}+32r^{3}+172r^{2} + 352r
	+ 15\left ( -1 \right )^{r} +225 \right )}{960} + \frac{2(r^{3}+12r^{2}+39r+28)}{54}  \Bigg] \\  + \frac{4r^{3}+30r^{2} + 68r -3 + 3 ( -1 )^{r}}{24} \text{ if }\, r  \bmod\!\!
\text{ } 3 = 2. \!\!\!\! \end{matrix}\right. \phantom{.......}
\end{eqnarray}

\begin{figure}[t!]
\centering{
\vspace{-5pt}
\includegraphics[scale=0.24]{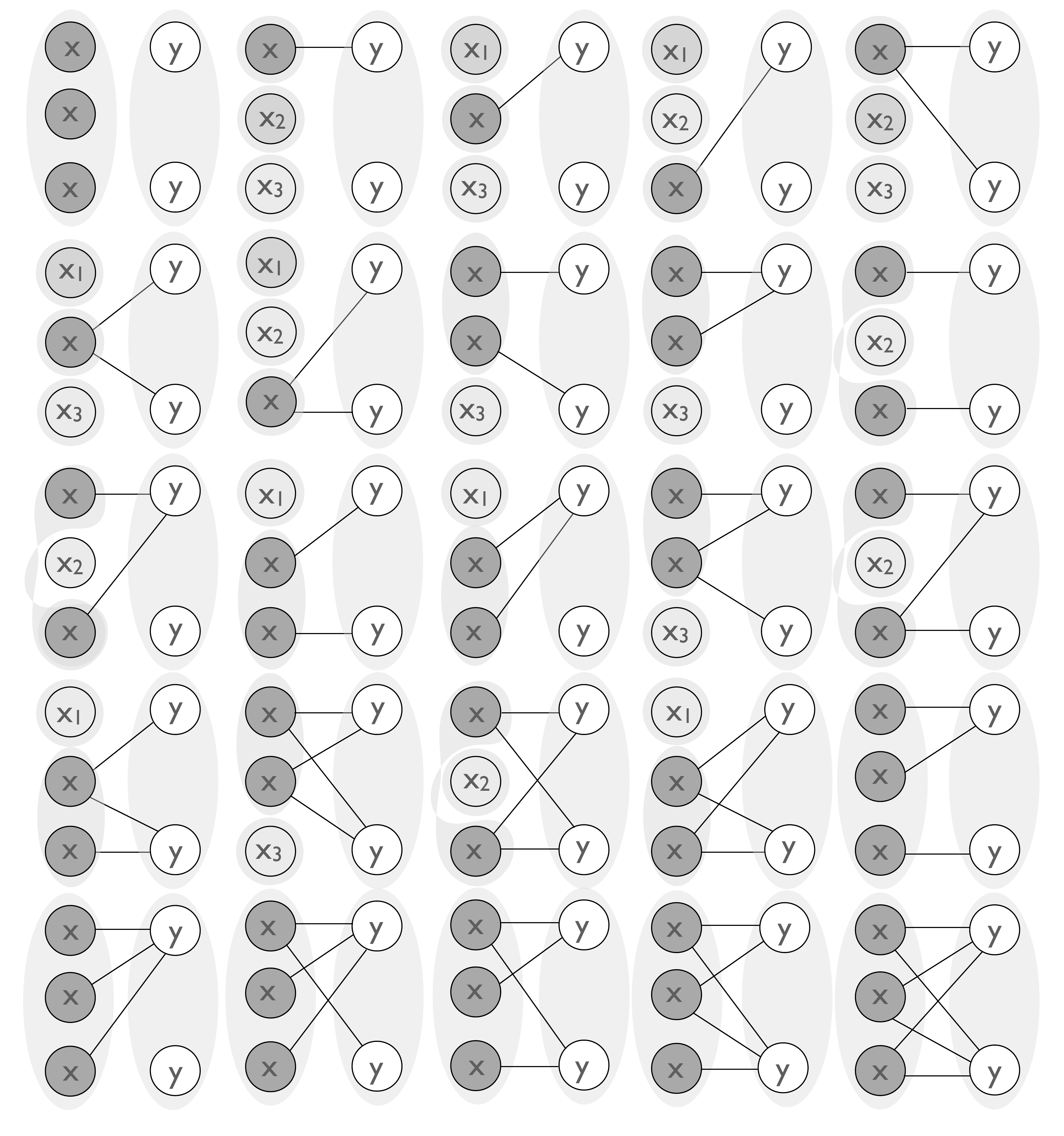}
\vspace{-10pt}
\caption{Non-equivalent left set-labeled (3,2)-bipartite graphs.}
\label{figure4}
}
\end{figure}

\vspace{-1pt}\noindent
The number of $(3,2)$-left-set-labeled bipartite graphs is determined to be 25 by substituting 3 for $r$  into the third formula above for $|B_x(3,r)|$ and all such graphs are depicted in Figure~\ref{figure4}.

Our counting method can also be used when the number of left vertices  $n$ is arbitrarily large, while the number of right vertices $r$ is  $2$ or $3$ as follows.
\vspace{-1pt}
\begin{eqnarray}
|\overline{B}_x(i,1)| &=& |B_u(i,1)| - |B_u(i-1,1)|, \\ 
	&=& (i+1) - i = 1.
\end{eqnarray}

\vspace{-5pt}\noindent 
This leads to

\vspace{-12pt}
\begin{equation}
|B_x(n,1)| = 1 + \sum_{i=1}^n \binom{n}{i} |\overline{B}_x(i,1)| = 2^n. 
\end{equation}

\vspace{2pt}\noindent
Similarly, $|\overline{B}_x(i,2)|$ can be calculated using (\ref{bxnrRecFormula}) as follows.
\vspace{-4pt}
\begin{eqnarray}
|\overline{B}_x(i,2)| &=& |B_u(i,2)| - |B_u(i-1,2)|, \\ 
&=& |B_u(2,i)| - |B_u(2, i-1)|, \\ 
&=& \frac{2i^{3}+15i^{2} + 34i + 22.5 + 1.5(-1)^{i}}{24}  \nonumber \\ 
&&\!\!\!\!\!- \frac{2(i-1)^{3}+15(i-1)^{2} + 34(i-1) + 22.5 + 1.5(-1)^{i-1}}{24}, \\
&=& \frac{6i^{2} + 24i + 21 + 3(-1)^{i}}{24}.
\label{bbarxn2}
\end{eqnarray}
\noindent This leads to

\vspace{-20pt}
\begin{eqnarray}
|B_x(n,2)| &=& 1 + \sum_{i=1}^n \binom{n}{i} |\overline{B}_x(i,2)|, \\
&=& 1 + \sum_{i=1}^n \binom{n}{i} \frac{6i^{2} + 24i + 21 + 3(-1)^{i}}{24}, \\
&=& n (n+1)2^{n-4}  + n 2^{n-1} + 7 (2^{n -3}).
\end{eqnarray}

\noindent
$|\overline{B}_x(i,3)|$ can similarly be computed as follows.

\begin{eqnarray}
|\overline{B}_x(i,3)| &=& |B_u(i,3)| - |B_u(i-1,3)| \\\\
&=& |B_u(3,i)| - |B_u(3, i-1)|  \\ \\
&=& \left\{\begin{matrix}
\frac{1}{6}\Bigg[ \binom{i+7}{7} - \binom{i+6}{7}  + \frac{3( i+4)( 2i^{4}+32i^{3}+172i^{2} + 352i + 15(-1)^{i} +225)}{960} \\
- \frac{3( i+3)( 2(i-1)^{4}+32(i-1)^{3}+172(i-1)^{2} + 352(i-1) - 15(-1)^{i} +225)}{960} \\
+ \frac{2(i^{3}+12i^{2}+45i+54)}{54} - \frac{2((i-1)^{3}+12(i-1)^{2}+39(i-1)+28)}{54} \Bigg] \text{ if }\, i \bmod\!\! \text{ } 3 = 0, \\
\frac{1}{6}\Bigg[ \binom{i+7}{7} - \binom{i+6}{7}  + \frac{3\left ( i+4 \right )\left ( 2i^{4}+32i^{3}+172i^{2} + 352i + 15 (-1 )^{i} +225 \right )}{960} \\
- \frac{3( i+3)( 2(i-1)^{4}+32(i-1)^{3}+172(i-1)^{2} + 352(i-1) - 15(-1)^{i} +225)}{960} \\
+ \frac{2(i^{3}+12i^{2}+45i+50)}{54} - \frac{2((i-1)^{3}+12(i-1)^{2}+45(i-1)+54)}{54} \Bigg] \text{ if }\, i  \bmod\!\! \text{ } 3 = 1, \\
\frac{1}{6}\Bigg[ \binom{i+7}{7} - \binom{i+6}{7}  + \frac{3(i+4)(2i^{4}+32i^{3}+172i^{2} + 352i + 15(-1)^{i} +225 )}{960} \\
- \frac{3( i+3)( 2(i-1)^{4}+32(i-1)^{3}+172(i-1)^{2} + 352(i-1) - 15(-1)^{i} +225)}{960} \\
+ \frac{2(i^{3}+12i^{2}+39i+28)}{54}  - \frac{2((i-1)^{3}+12(i-1)^{2}+45(i-1)+50)}{54} \Bigg] \text{ if }\, i  \bmod\!\!  \text{ } 3 = 2, \!\!\!\! \end{matrix}\right. \phantom{.......}
\end{eqnarray}
\vspace{-5pt}
\begin{eqnarray}
&=& \left\{\begin{matrix}
\frac{1}{6}\Bigg[ \binom{i+6}{6}  + \frac{10i^{4}+140i^{3}+680i^{2} + 1330i + 30i(-1)^{i}+ 105(-1)^{i} +855}{320} \\
+ \frac{6i^{2}+54i+108}{54} \Bigg] \text{ if }\, i \bmod\!\! \text{ } 3 = 0, \\
\frac{1}{6}\Bigg[ \binom{i+6}{6} + \frac{10i^{4}+140i^{3}+680i^{2} + 1330i + 30i(-1)^{i}+ 105(-1)^{i} +855}{320} \\
+ \frac{6i^{2}+42i+60}{54} \Bigg] \text{ if }\, i  \bmod\!\! \text{ } 3 = 1, \\
\frac{1}{6}\Bigg[ \binom{i+6}{6} + \frac{10i^{4}+140i^{3}+680i^{2} + 1330i + 30i(-1)^{i}+ 105(-1)^{i} +855}{320} \\
+ \frac{6i^{2}+30i+24}{54} \Bigg] \text{ if }\, i  \bmod\!\!  \text{ } 3 = 2. \!\!\!\! \end{matrix}\right.
\label{bbarxn3}
\end{eqnarray}

\noindent This leads to

\vspace{-20pt}
\begin{eqnarray}
|B_x(n,3)| &=& 1 + \sum_{i=1}^n \binom{n}{i} |\overline{B}_x(i,3)|, \\\\x
&=& \left\{\begin{matrix}
\frac{3\ 2^n n^6+171\ 2^n n^5+3765\ 2^n n^4+41265\ 2^n n^3+14787\ 2^{n+4} n^2}{829440} \\
+ \frac{12 \left(2560 (-1)^{n/3}+55077\ 2^n\right) n+880 \left(128 (-1)^{n/3}+763\
	2^n\right)}{829440} \text{ if }\, n \bmod\!\! \text{ } 3 = 0, \\
\\
\frac{3\ 2^n n^6+171\ 2^n n^5+3765\ 2^n n^4+41265\ 2^n n^3+14787\ 2^{n+4} n^2}{829440} \\
+ \frac{165231\ 2^{n+2} n-80 \left(1280 (-1)^{\frac{n+2}{3}}-8393\ 2^n\right)}{829440} \text{ if }\, n \bmod\!\! \text{ } 3 = 1, \\
\\
\frac{3\ 2^n n^6+171\ 2^n n^5+3765\ 2^n n^4+41265\ 2^n n^3+14787\ 2^{n+4} n^2}{829440} \\
+\frac{12 \left(2560 (-1)^{\frac{n+1}{3}}+55077\ 2^n\right) n+80 \left(128
	(-1)^{\frac{n+1}{3}}+8393\ 2^n\right)}{829440} \text{ if }\, n \bmod\!\!  \text{ } 3 = 2, \!\!\!\! \end{matrix}\right.
\label{bxn3Sum}
\end{eqnarray}

\vspace{-5pt}\noindent where $n \ge 2$.

\vspace{2pt}
\section{Set-Labeled Bipartite Graphs}
\label{setLabel}

\vspace{-6pt}\noindent
We now count  set-labeled $(n,2)$ and $(n,3)$-bipartite graphs for $n\ge 1.$
Let $\overline{B}_{x,y}(i,j)$ be the set of all $(i,j)$-unlabeled bipartite graphs such that there is no vertex in the graph that has a degree of 0. We have\footnote{Here, we could use a notation that is similar to the notation $B_{x}(n,r,i)$ of the prior section for completeness. However, the relation between $\overline{B}_{x,y}(i,j)$ and  $B_{x,y}(n,r,i,j)$ is obvious and omitted here.}
\begin{equation}
|B_{x,y}(n,r)| = 1 + \sum_{i=1}^n \sum_{j=1}^r \binom{n}{i} \binom{r}{j} |\overline{B}_{x,y}(i,j)|.
\label{bx,ySumFormula}
\end{equation}

\vspace{-5pt}\noindent
Clearly, $|\overline{B}_{x,y}(1,1)| = |\overline{B}_{x,y}(1,j)| = |\overline{B}_{x,y}(i,1)| = 1$. $|\overline{B}_x(i,j)$ and $|\overline{B}_x(i,j-1)|$ can be witten as follows.
\vspace{-15pt}
\begin{eqnarray}
|\overline{B}_x(i,j)| = \sum_{k=1}^j |\overline{B}_{x,y}(i,k)|,  \\
|\overline{B}_x(i,j-1)| = \sum_{k=1}^{j-1} |\overline{B}_{x,y}(i,k)|.
\end{eqnarray}

\vspace{-10pt}\noindent
Hence
\begin{equation}
|\overline{B}_{x,y}(i,j)| = |\overline{B}_x(i,j)| - |\overline{B}_x(i,j-1)|.
\label{bx,yRecursion}
\end{equation}

\noindent\vspace{5pt}

Now we can calculate $|\overline{B}_{x,y}(i,2)|$ using Eqn.~\ref{bx,yRecursion}.

\begin{eqnarray}
|\overline{B}_{x,y}(i,2)| & = & |\overline{B}_x(i,2)| - |\overline{B}_x(i,1)|, \\\\
&= & \frac{6i^{2} + 24i + 21 + 3(-1)^{i}}{24} - 1, \\\\
&= & \frac{6i^{2} + 24i -3 + 3(-1)^{i}}{24}.
\label{bbarx,yn2}
\end{eqnarray}

\vspace{5pt}\noindent 
Using Eqns.~\ref{bx,ySumFormula} and \ref{bbarx,yn2}, we can calculate $|B_{x,y}(n,2)|$ as follows.

\vspace{-15pt}
\begin{eqnarray}
|B_{x,y}(n,2)| &=& 1 + \sum_{i=1}^n \binom{n}{i} \sum_{j=1}^2 \binom{2}{j} |\overline{B}_{x,y}(i,j)|, \\\\
&=& 1 + \sum_{i=1}^n \binom{n}{i} \Bigg[ \binom{2}{1} |\overline{B}_{x,y}(i,1)| + \binom{2}{2} |\overline{B}_{x,y}(i,2)|  \Bigg], \\\\
&=& 1 + \sum_{i=1}^n \binom{n}{i} \frac{6i^{2} + 24i +45 + 3(-1)^{i}}{24}, \\
&=& \frac{15}{8} \left(2^n\right)+2^{n-1} n+2^{n-4} n (n+1)-1.
\end{eqnarray}

\begin{figure}[t!]
\centering{
\vspace{5pt}
\includegraphics[scale=0.22]{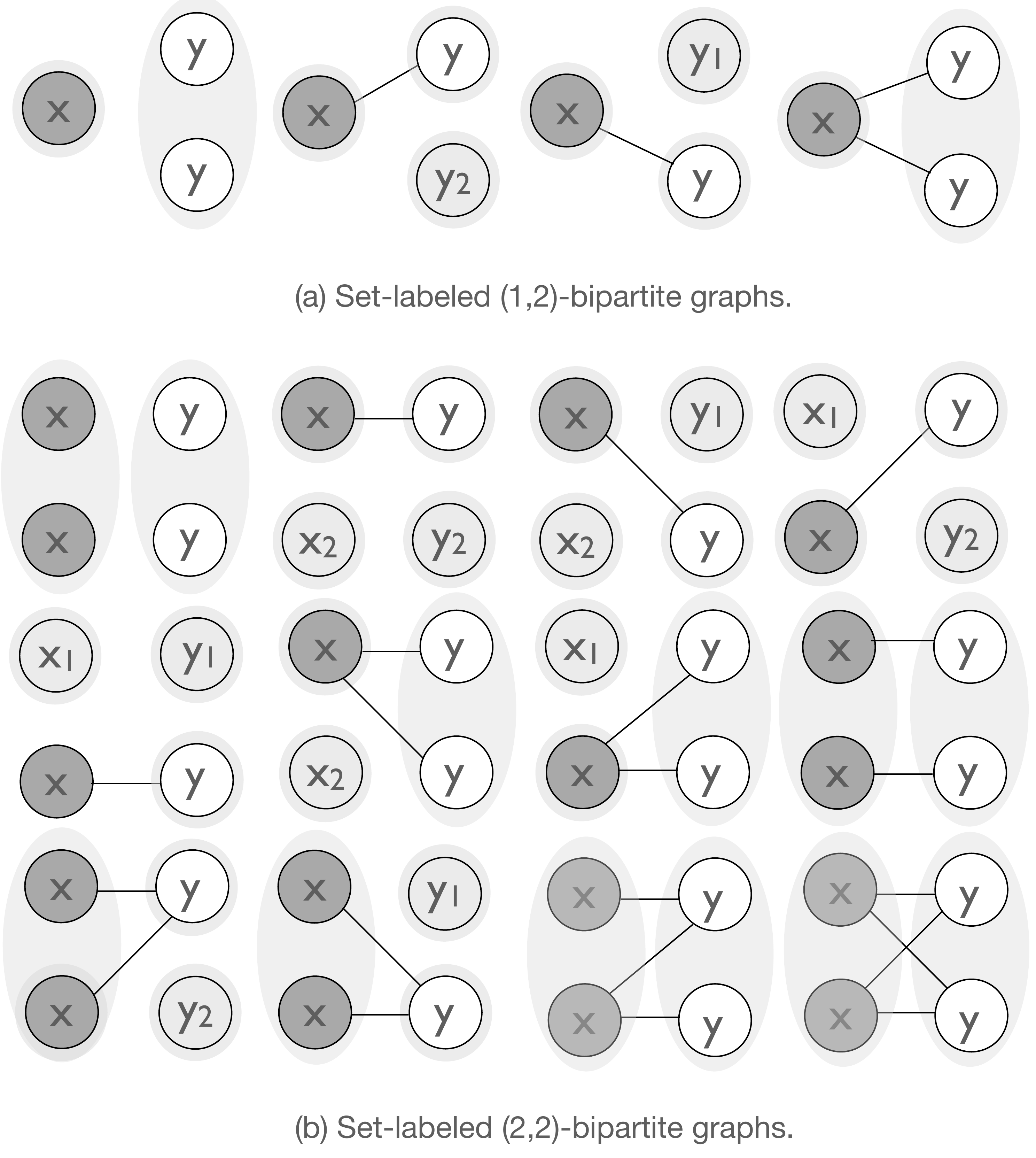}
\vspace{-10pt}
\caption{Non-equivalent set-labeled (1,2) and (2,2)-bipartite graphs.}
\label{figure5}
}
\end{figure}

\begin{figure}[t!]
\centering{
\vspace{-5pt}
\includegraphics[scale=0.19]{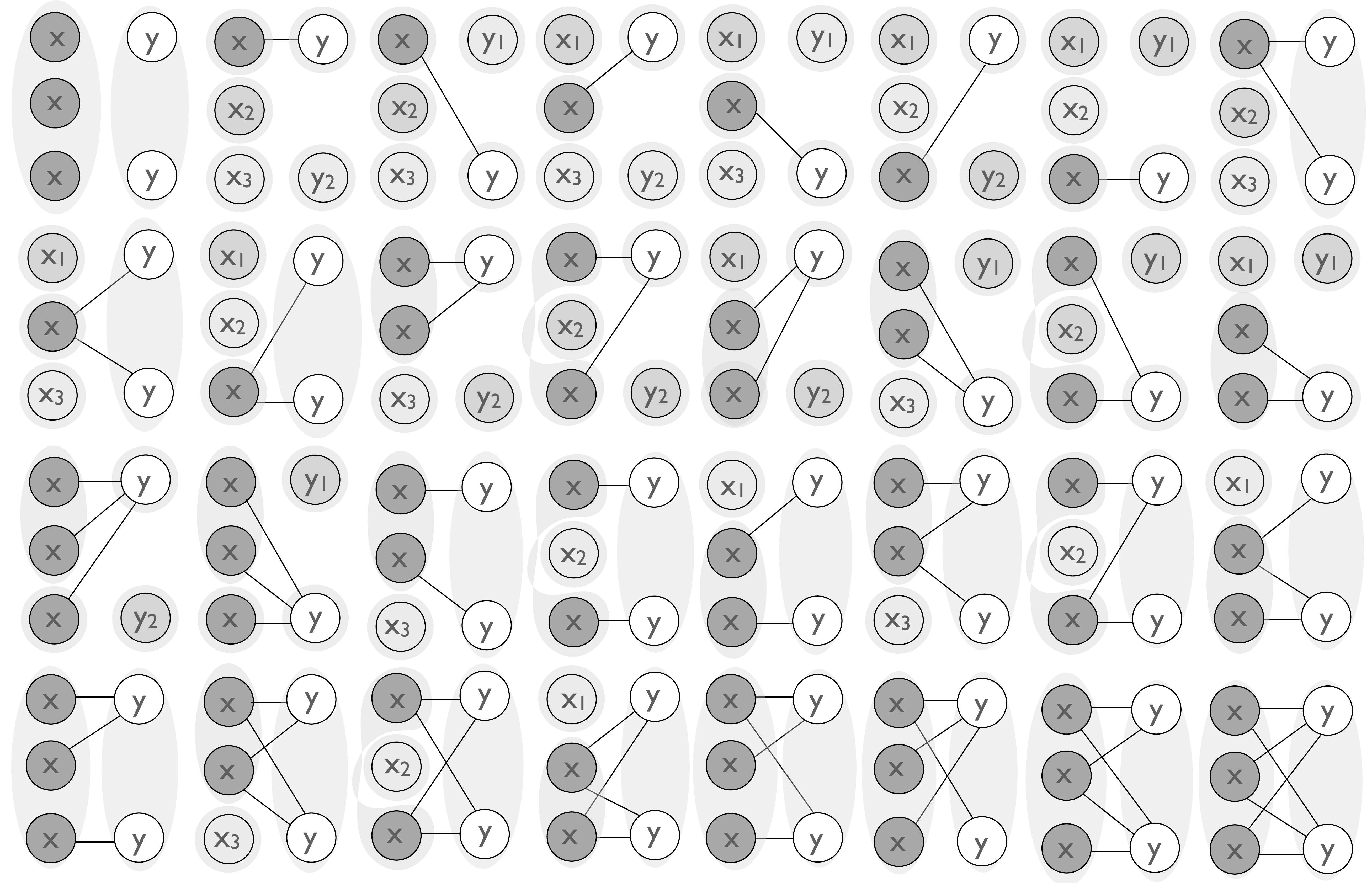}
\vspace{-10pt}
\caption{Non-equivalent set-labeled (3,2)-bipartite graphs.}
\label{figure6}
}
\end{figure}

\begin{remark}
{\rm
It follows that $|B_{x,y}(1,2)| = 4,$ $|B_{x,y}(2,2)| = 12$ and $|B_{x,y}(3,2)| = 32.$ The corresponding bipartite graphs in $B_{x,y}(1,2), B_{x,y}(2,2)$ and  $B_{x,y}(3,2)$ are shown in Figure~\ref{figure5} and Figure~\ref{figure6}. $||$
}
\end{remark}

\vspace{5pt}\noindent
Similarly, we can calculate $|\overline{B}_{x,y}(i,3)|$ by letting $j = 3$ in Eqn.~\ref{bx,yRecursion} to get 

\vspace{5pt}
\begin{equation}
|\overline{B}_{x,y}(i,3)| = |\overline{B}_x(i,3)| - |\overline{B}_x(i,2)|.
\label{bbarx,yi3}
\end{equation}

\vspace{5pt}\noindent 
Now we use Eqns.~\ref{bbarxn2} and \ref{bbarx,yi3} to compute $|B_{x,y}(n,3)|$ as follows.

\vspace{-10pt}
\begin{eqnarray}
|B_{x,y}(n,3)| &=& 1 + \sum_{i=1}^n \binom{n}{i} \sum_{j=1}^3 \binom{3}{j} |\overline{B}_{x,y}(i,j)|, \\
&=& 1 + \sum_{i=1}^n \binom{n}{i} \Bigg[ \binom{3}{1} |\overline{B}_{x,y}(i,1)| + \binom{3}{2} |\overline{B}_{x,y}(i,2)| + \binom{3}{3} |\overline{B}_{x,y}(i,3)|  \Bigg], \\
&=& 1 + \sum_{i=1}^n \binom{n}{i} \Bigg[ 3 + 3 |\overline{B}_{x,y}(i,2)| +  |\overline{B}_x(i,3)| - |\overline{B}_x(i,2)|  \Bigg], \\
&=& 1 + \sum_{i=1}^n \binom{n}{i} \Bigg[ 3 + \frac{3(6i^{2} + 24i -3 + 3(-1)^{i})}{24} - \frac{6i^{2} + 24i + 21 + 3(-1)^{i}}{24}  \Bigg] \nonumber \\
& \text{ }& \,\,\, + \sum_{i=1}^n \binom{n}{i} |\overline{B}_x(i,3)|, \\
& = &1 + \sum_{i=1}^n \binom{n}{i} \Bigg[ \frac{i^2}{2}+2 i+\frac{(-1)^i}{4}+\frac{7}{4}  \Bigg] +  \sum_{i=1}^n \binom{n}{i} |\overline{B}_x(i,3)|, \\
& = & 1+ 2^{n-3} n^2+9\times 2^{n-3} n+7\times 2^{n-2}-2 +  \sum_{i=1}^n \binom{n}{i} |\overline{B}_x(i,3)|.
\end{eqnarray}

\vspace{-2pt}\noindent 
Rearranging terms and using Eqn.~\ref{bxn3Sum} gives

\vspace{-10pt}\begin{eqnarray}
|B_{x,y}(n,3)| &=& 2^{n-3} n^2+9\times 2^{n-3} n+7\times 2^{n-2}-2 + |B_{x}(n,3)| \\ \\
&=& \left\{\begin{matrix}
\frac{3\times 2^n n^6+171\times 2^n n^5+3765\times 2^n n^4+41265\times 2^n n^3+21267\times 2^{n+4} n^2}{829440}\\
+  \frac{12 \left(2560 (-1)^{n/3}+132837\times 2^n\right) n+80
	\left(1408 (-1)^{n/3}+26537\times 2^n-20736\right)}{829440} \text{ if }\, n \bmod\!\! \text{ } 3 = 0, \\
\\
\frac{3\times 2^n n^6+171\times 2^n n^5+3765\times 2^n n^4+41265\times 2^n n^3+21267\times 2^{n+4} n^2+398511\times 2^{n+2} n}{829440} \\
- \frac{80 \left(1280
	(-1)^{\frac{n+2}{3}}-26537\times 2^n+20736\right)}{829440} \text{ if }\, n \bmod\!\! \text{ } 3 = 1, \\
\\
\frac{3\times 2^n n^6+171\times 2^n n^5+3765\times 2^n n^4+41265\times 2^n n^3+21267\times 2^{n+4} n^2}{829440} \\
+ \frac{12 \left(2560 (-1)^{\frac{n+1}{3}}+132837\times
	2^n\right) n+80 \left(128 (-1)^{\frac{n+1}{3}}+26537\times 2^n-20736\right)}{829440}  \text{ if }\, n \bmod\!\!  \text{ } 3 = 2. \!\!\!\! \end{matrix}\right.
\end{eqnarray}

\vspace{-10pt}
\section{Lower And Upper Bounds for $|B_x(n,r)|$ }

\vspace{-10pt}

We now present lower and upper bounds for the number of left-set-labeled bipartite graphs for all  $n, r\ge 3.$ 

\vspace{2pt}
\begin{theorem}
$\,$\\
\vspace{-20pt}
	\label{propLowBoundSet}
	{\rm	
		\begin{eqnarray}
		|B_x(n,r)|  \ge 1 + \sum_{i=1}^n \binom{n}{i} \binom{i+2^r-2}{i}\text{{\Large/}}r!.
		\end{eqnarray}
	}
\end{theorem}

\vspace{-5pt}
\begin{mproof}
	{\rm
		The first term $1$ given in lower bound formula counts the left-set-labeled bipartite graph without any edges. To count the remaining graphs in $B_x(n,r),$ consider $2^r - 1$ subsets of $r$ right vertices excluding the empty set. Let $X_i$ be an arbitrary but fixed subset of left $n$ vertices of size $i, 1\le i \le n$.  Each one-to-one pairing of $X_i$ with any $i$ of the $2^r -1 $ subsets of right vertices constitutes a left-set-labeled bipartite graph. The number of such pairings is given by\footnote{This formula counts the number $i$-selections from a set of $2^r-1$ distinct elements.} $\binom{i+2^r-2}{i}$. Given that right vertices are indistinguishable, possibly up to $r!$, but no more of these left-set-labeled graphs may be equivalent under a permutation of $r$ vertices. Therefore, there exist at least $\binom{i+2^r-2}{i}\text{{\large/}}r!$ distinct left-set-label bipartite graphs associated with $X_i$. Since there exist $\binom{n}{i}$ $i$-subsets of left vertices, we must have at least $\binom{n}{i}\binom{i+2^r-2}{i}\text{{\large/}}r!$ distinct left-set-labeled bipartite graphs between subsets of  $i$ left vertices and the $r$ right vertices. Summing this for $i$ from $1$ to $n$ establishes the lower bound in the statement.\qed
	}
\end{mproof}

It is difficult to obtain a closed form formula for the sum in Theorem~\ref{propLowBoundSet}, but it can be approximated by computing the maximum value of the argument of the sum with respect to $i$ as shown in the following proposition.

\vspace{5pt}
\begin{corollary}
	{\rm
	$\,$\\
	\vspace{-20pt}
		\begin{eqnarray}
		|B_x(n,r)|  \ge \binom{n}{i_{max}} \binom{i_{max}+2^r-2}{i_{max}}\text{{\Large/}}r!,
		\end{eqnarray}
		
		\vspace{-1pt}\noindent
		where $i_{max} = \frac{1}{4} \Big[ (a^2 + 6an + n^2)^{1/2} -a-n \Big]$ and $a= 2^r - 2$.
	}
\end{corollary}
\vspace{8pt}
\begin{mproof}
	{\rm
The inequality
\vspace{-2pt}
\begin{equation}
	|B_x(n,r)|  \ge  \binom{n}{i_{max}} \binom{i_{max}+2^r-2}{i_{max}}\text{{\Large/}}r! 
	\end{equation}
	
	\vspace{4pt}\noindent
	holds, where $i_{max}$ is the value of $i$, $1 \le i_{max} \le n$ that maximizes the expression $\binom{n}{i} \binom{i+2^r-2}{i}\text{{\Large/}}r!$.
	
	\noindent
	Differentiating the argument of the sum and setting it equal to 0 gives 
	$$ i_{max} = \frac{1}{4} \Big[ (a^2 + 6an + n^2)^{1/2} -a-n \Big] + n/2,$$ 
	
	\vspace{-3pt}\noindent
	where $a= 2^r - 2$. \qed
	}
\end{mproof}

\vspace{5pt}
The next theorem provides an upper bound on $|B_x(n,r)|.$

\vspace{5pt}
\begin{theorem}
	{\rm
	$\,$\\
	\vspace{-25pt}
\begin{eqnarray}
 |B_x(n,r)| \le 1+\sum_{i=1}^n \binom{n}{i} \Bigg[ 2 \binom{r+2^i-1}{r}\text{{\Large/}}i! - \binom{r+2^{i-1}-1}{r}\text{{\Large/}}(i-1)! \Bigg]\!,\text{ }
\end{eqnarray}

\vspace{-10pt}\noindent
where $n < r $.
}
\end{theorem}

\vspace{4pt}
\begin{mproof}
	{\rm
	As before, the first term counts the left-set-labeled bipartite graph with no edges. For the remaining sum, we use Eqns.~\ref{bxnrFormula}  and \ref{bxnrRecFormula}, and the lower and upper bounds in Eqn.~\ref{generalCase} to obtain
	
	\begin{eqnarray}
	|B_x(n,r)| &\!\!\!\!= &\!\!\!\! 1+ \sum_{i=1}^{n} \binom{n}{i} |\overline{B}_x(i,r)|,  \\ 
	&\!\!\!\! =&\!\!\!\! 1+ \sum_{i=1}^{n} \binom{n}{i} \Big[ |B_u(i,r)| - |B_u(i-1,r)| \Big], \\
	&\!\!\!\! \le &\!\!\!\! 1+\sum_{i=1}^n \binom{n}{i} \Bigg[ 2 \binom{r+2^i-1}{r}\text{{\Large/}}i! - \binom{r+2^{i-1}-1}{r}\text{{\Large/}}(i-1)! \Bigg], \phantom{........}
	\end{eqnarray}
	and this establishes the upper bound in the statement.\qed
}
\end{mproof}

\vspace{5pt}
\begin{theorem}
	{\rm For all $n, r \ge 2$,
	
	\vspace{-20pt}
	\begin{eqnarray}
	|B_x(n,r)| & \le & 1+2 \sum_{i=1}^{r-1} \binom{n}{i} \binom{r+2^i-1}{r} \text{{\Large/}}i!+ 2 \binom{n}{r}  \binom{r+2^r}{r+1}\text{{\Large/}}r! \nonumber \\
	& & +2 \sum_{i=r+1}^n \binom{n}{i} \binom{i+2^r-1}{i}\text{{\Large/}}r!.
	\end{eqnarray}	
	}
\end{theorem}

\begin{mproof}
	{\rm
		Again the first term $1$ given in the upper bound formula counts for the left-set-labeled bipartite graph where there are no edges between vertices. For the remaining sum, we note $|\overline{B}_{x}(i,r)|$ can not exceed $|B_u(i,r)|$, since every bipartite graph in $\overline{B}_{x}(i,r)$ is also in $B_u(i,r)$, $1\le i \le n $.
		\begin{eqnarray}
		|B_x(n,r)| &= & 1+ \sum_{i=1}^{n} \binom{n}{i} |\overline{B}_x(i,r)|, \\ 
		& \le & 1+ \sum_{i=1}^{n} \binom{n}{i}  |B_u(i,r)|, \\
		& \le & 1+2 \sum_{i=1}^{r-1} \binom{n}{i} \binom{r+2^i-1}{r}\text{{\Large/}}i!+ \sum_{i=r}^{r}\binom{n}{i} |B_u(i,r)| \nonumber \\
		 & & +2 \sum_{i=r+1}^n \binom{n}{i} \binom{i+2^r-1}{i}\text{{\Large/}}r!.
		\end{eqnarray}
		
		\noindent Since we do not have any upper bound formula for $|B_u(i,r)|$ when $i=r$, we will use $|B_u(r,r)| \le |B_u(r+1,r)|$.
		\begin{eqnarray}
		|B_x(n,r)| & \le & 1+2 \sum_{i=1}^{r-1} \binom{n}{i} \binom{r+2^i-1}{r} \text{{\Large/}}i!+ 2 \binom{n}{r}  \binom{r+2^r}{r+1}\text{{\Large/}}r! \nonumber \\
		& & +2 \sum_{i=r+1}^n \binom{n}{i} \binom{i+2^r-1}{i}\text{{\Large/}}r!,
		\end{eqnarray}	
		\noindent
		and this establishes the upper bound in the statement.\qed
	}
\end{mproof}

\section{Lower and Upper Bounds For  $|B_{x,y}(n,r)|$}

In this section, we provide both a lower and upper bound on the number of set-labeled bipartite graphs.  
\noindent
We first use  the 2-sided inequality\cite{atmacaoruc2018}

\begin{equation}
\label{atmacaOrucInequality}
\displaystyle  \frac{\binom{r+2^{n}-1}{r}}{n!}\, \le\, |B_u(n,r)|\, \le\,  2\frac{\binom{r+2^{n}-1}{r}}{n!},\, n <  r.\nonumber
\end{equation}
 
and the lower bound (See Section 4, p. 712 in\cite{atmacaoruc2018})

\begin{equation}
\label{atmacaOrucInequality2}
\displaystyle  \frac{\binom{n+2^{n}-1}{n}}{2n!}\, \le\, |B_u(n,n)|\nonumber
\end{equation}

to obtain following inequalities for $|B_u(i,j)|$ and $|B_u(i-1,j)|.$
\begin{eqnarray}
\frac{ \binom{j+2^{i}-1}{j}}{i!} \le & |B_u(i,j)| & \le  \frac{2\binom{j+2^{i}-1}{j}}{i!}, \\
\frac{ \binom{j+2^{i-1}-1}{j}}{(i-1)!}  \le & |B_u(i-1,j)| & \le  \frac{2\binom{j+2^{i-1}-1}{j}}{(i-1)!},
\end{eqnarray}
\noindent where $i<j$,
	\begin{eqnarray}
	\frac{ \binom{i+2^{j}-1}{i}}{j!} \le & |B_u(i,j)| & \le  \frac{2\binom{i+2^{j}-1}{i}}{j!}, \\
	\frac{ \binom{i+2^{j}-2}{i-1}}{j!} \le & |B_u(i-1,j)| & \le  \frac{2\binom{i+2^{j}-2}{i-1}}{j!},
	\end{eqnarray}
\noindent where $j<i-1$,
\begin{eqnarray}
\frac{ \binom{i+2^{i-1}-1}{i}}{(i-1)!} \le & |B_u(i,j)| & \le  \frac{2\binom{i+2^{i-1}-1}{i}}{(i-1)!},  \\
\frac{ \binom{i+2^{i-1}-2}{i-1}}{2(i-1)!} \le & |B_u(i-1,j)| & \le  \frac{2\binom{i+2^{i-1}-1}{i}}{(i-1)!}, 
\end{eqnarray}
\noindent where $j=i-1$ and
\begin{eqnarray}
\frac{ \binom{i+2^{i}-1}{i}}{2(i!)} \le & |B_u(i,j)| & \le  \frac{2\binom{i+2^{i}}{i+1}}{i!},  \\
\frac{ \binom{i+2^{i-1}-1}{i}}{(i-1)!} \le & |B_u(i-1,j)| & \le  \frac{2\binom{i+2^{i-1}-1}{i}}{(i-1)!}, 
\end{eqnarray}
\noindent where $j=i$.

\noindent
Using these inequalities with Eqn.~\ref{bxnrRecFormula}, we can bound $|\overline{B}_x(i,j)|$ as follows.
\begin{eqnarray}
	\frac{\binom{j+2^{i}-1}{j}}{i!} - \frac{ 2\binom{j+2^{i-1}-1}{j}}{(i-1)!} & \le 	|\overline{B}_x(i,j)| & \le \frac{2\binom{j+2^{i}-1}{j}}{i!} - \frac{ \binom{j+2^{i-1}-1}{j}}{(i-1)!}, i<j, \\\\
	\frac{ \binom{i+2^{j}-1}{i}}{j!}- \frac{2\binom{i+2^{j}-2}{i-1}}{j!} & \le 	|\overline{B}_x(i,j)| & \le  \frac{2\binom{i+2^{j}-1}{i}}{j!} - \frac{ \binom{i+2^{j}-2}{i-1}}{j!}, j<i-1,\\\\
	\frac{ \binom{i+2^{i-2}-1}{i}}{(i-2)!}- \frac{2\binom{i+2^{i-2}-2}{i-1}}{(i-2)!} \!\!\!\!\!\!\!\! & \le |\overline{B}_x(i,i-2)| \le |\overline{B}_x(i,j)| \!\!\!\!\!\!&\!\! \le \!\!\! \frac{2\binom{i+2^{i-1}-1}{i}}{(i-1)!} - \frac{ \binom{i+2^{i-1}-2}{i-1}}{2(i-1)!}, j=i-1, \phantom{...........}\\\\
	\frac{\binom{i+2^{i}-2}{i-1}}{i!} - \frac{ 2\binom{i+2^{i-1}-2}{i-1}}{(i-1)!} & \le |\overline{B}_x(i,i-1)| \le |\overline{B}_x(i,j)| & \le \frac{2\binom{i+2^{i}}{i+1}}{i!} - \frac{ \binom{i+2^{i-1}-1}{i}}{(i-1)!}, j = i.
\end{eqnarray}

\begin{proposition}
	{\rm
		Let $B_{x,y}(n,r)$ be the set of all $(n,r)$-set-labeled bipartite graphs. 
		
		\begin{eqnarray}
		|B_{x,y}(n,r)| \!\!\!\!&\ge\!\!\!\!& 1+  \sum_{i=1}^n \binom{n}{i} \Bigg( 
		\sum_{j=1}^{i-2} \binom{r}{j} \Bigg[ \frac{ \binom{i+2^{j}-1}{i}}{j!}- \frac{2\binom{i+2^{j}-2}{i-1}}{j!} - \frac{2\binom{i+2^{j-1}-1}{i}}{(j-1)!} + \frac{ \binom{i+2^{j-1}-2}{i-1}}{(j-1)!}  \Bigg] \nonumber \\ 
		& \text{ } & +\sum_{j=i+2}^r \binom{r}{j} \Bigg[ \frac{\binom{j+2^{i}-1}{j}}{i!} - \frac{ 2\binom{j+2^{i-1}-1}{j}}{(i-1)!} - \frac{2\binom{j+2^{i}-2}{j-1}}{i!} + \frac{ \binom{j+2^{i-1}-2}{j-1}}{(i-1)!} \Bigg]
		\Bigg).
		\end{eqnarray}
	}
\end{proposition}

\begin{mproof}
	{\rm
	Substituting $|\overline{B}_{x,y}(i,j)|$ in Eqn.~\ref{bx,ySumFormula} using Eqn.~\ref{bx,yRecursion} gives
	
	\vspace{-8pt}
	\begin{eqnarray}
	|B_{x,y}(n,r)| \!\!\!\!&= \!\!\!\!& 1+ \sum_{i=1}^n \sum_{j=1}^r \binom{n}{i} \binom{r}{j} |\overline{B}_{x,y}(i,j)|,  \\\\
	& = & 1+  \sum_{i=1}^n \sum_{j=1}^r \binom{n}{i} \binom{r}{j} \Bigg[  |\overline{B}_x(i,j)| - |\overline{B}_x(i,j-1)| \Bigg].  
	\end{eqnarray}
	
	\vspace{-2pt}\noindent 
	Now replacing $|\overline{B}_x(i,j)|$ with its lower bound, $|\overline{B}_x(i,j-1)|$ with its upper bound and ignoring the cases where $i-1 \le j \le i+1$ gives
	
	\vspace{-5pt}
	\begin{eqnarray}
	|B_{x,y}(n,r)| \!\!\!\!&\ge\!\!\!\!& 1+  \sum_{i=1}^n \binom{n}{i} \Bigg( 
	\sum_{j=1}^{i-2} \binom{r}{j} \Bigg[ \frac{ \binom{i+2^{j}-1}{i}}{j!}- \frac{2\binom{i+2^{j}-2}{i-1}}{j!} - \frac{2\binom{i+2^{j-1}-1}{i}}{(j-1)!} + \frac{ \binom{i+2^{j-1}-2}{i-1}}{(j-1)!}  \Bigg] \nonumber \\ 
	& \text{ } & +\sum_{j=i+2}^r \binom{r}{j} \Bigg[ \frac{\binom{j+2^{i}-1}{j}}{i!} - \frac{ 2\binom{j+2^{i-1}-1}{j}}{(i-1)!} - \frac{2\binom{j+2^{i}-2}{j-1}}{i!} + \frac{ \binom{j+2^{i-1}-2}{j-1}}{(i-1)!} \Bigg]	\Bigg).
	\end{eqnarray}
\qed}
\end{mproof}

\vspace{-5pt}
\begin{remark}
	\rm{
	\noindent 
	It is noted that using the cases where $i-1 \le j \le i+1$ will likely give negative terms. Therefore, they are not included in the computation of the lower bound.
}
\end{remark}

\vspace{5pt}
\begin{proposition}
		{\rm
			Let $B_{x,y}(n,r)$ be the set of all $(n,r)$-set-labeled bipartite graphs. 
			
			\vspace{-20pt}
			\begin{eqnarray}
			|B_{x,y}(n,r)| \!\!\!\! & \le \!\!\!\!& 1+2\sum_{i=1}^n \binom{n}{i} \Bigg[ \sum_{j=1}^{i-1} \binom{r}{j} \binom{i+2^j-1}{i}\text{{\Large/}}j! + \binom{r}{i} \binom{i+2^{i}}{i+1}\text{{\Large/}}i! \nonumber \\
			& \text{\,\,\,\,} &+ \sum_{j=i+1}^r \binom{r}{j} \binom{j+2^i-1}{j}\text{{\Large/}}i! \Bigg]. 
			\end{eqnarray}
		}
	\end{proposition}
	
	\vspace{-10pt}
	\begin{mproof}
		{\rm
	
\noindent For the upper bound, we note that $|\overline{B}_{x,y}(i,j)|$ can not exceed $|B_u(i,j)|$, since every bipartite graph in $\overline{B}_{x,y}(i,j)$ is also in $B_u(i,j)$, $1\le i \le n $, $1\le j \le r $. 
	
\vspace{-18pt}
\begin{eqnarray}
	|B_{x,y}(n,r)| \!\!\!\!&= \!\!\!\!& 1+ \sum_{i=1}^n \sum_{j=1}^r \binom{n}{i} \binom{r}{j} |\overline{B}_{x,y}(i,j)|,  \\ 
	\!\!\!\!& \le \!\!\!\!& 1+  \sum_{i=1}^n \sum_{j=1}^r \binom{n}{i} \binom{r}{j}  |B_u(i,j)|,  \\ 
	\!\!\!\!& \le \!\!\!\!& 1+2\sum_{i=1}^n \binom{n}{i} \Bigg[ \sum_{j=1}^{i-1} \binom{r}{j} \binom{i+2^j-1}{i}\text{{\Large/}}j! + \binom{r}{i} \binom{i+2^{i}}{i+1}\text{{\Large/}}i! \nonumber \\
	& \text{\,\,\,\,} &+ \sum_{j=i+1}^r \binom{r}{j} \binom{j+2^i-1}{j}\text{{\Large/}}i! \Bigg],
	\end{eqnarray}
	
	\vspace{-5pt}\noindent
	gives the formula in the statement. Note that the second term in the last inequality has been obtained by replacing $j=i$ by $j=i+1$.\qed
}
\end{mproof}

\bibliographystyle{unsrt}

\end{document}